\documentclass[11pt]{article}
\usepackage{latexsym}
\usepackage{amsfonts}
\usepackage{mathrsfs,amsthm,amsmath}
\usepackage{color}
\usepackage{todonotes}
\newtheorem{theorem}{Theorem}[section]
\newtheorem{proposition}{Proposition}[section]

\newtheorem{remark}{Remark}[section]

\begin{document}
\title{Noncentral moderate deviations for fractional Skellam processes}
\author{Jeonghwa Lee\thanks{Address: Department of Mathematics and Statistics,
University of North Carolina Wilmington, Sartarelli Hall, Office 2012F, 601 South College Road,
Wilmington, NC 28403-5970, USA. e-mail: \texttt{leejb@uncw.edu}} \and
Claudio Macci\thanks{Address: Dipartimento di Matematica, Universit\`a di Roma Tor 
Vergata, Via della Ricerca Scientifica, I-00133 Rome, Italy. e-mail:
\texttt{macci@mat.uniroma2.it}}}
\maketitle
\begin{abstract}
The term \emph{moderate deviations} is often used in the literature to mean a class of large deviation 
principles that, in some sense, fills the gap between a convergence in probability to zero (governed by 
a large deviation principle) and a weak convergence to a centered Normal distribution. We talk about 
\emph{noncentral moderate deviations} when the weak convergence is towards a non-Gaussian distribution.
In this paper, we present noncentral moderate deviation results for two fractional Skellam processes 
in the literature (see \cite{KerssLeonenkoSikorskii}). We also establish that, for the fractional Skellam
process of type 2 (for which we can refer to the recent results for compound fractional Poisson processes in
\cite{BeghinMacciSPL2022}), the convergences to zero are usually faster because we can prove suitable 
inequalities between rate functions.
\end{abstract}

\section{Introduction}
A large deviation principle provides some asymptotic bounds for a family of probability measures on the same topological space 
$\mathcal{X}$; moreover one often refers to a family of $\mathcal{X}$-valued random variables, $\{C_t\}$, whose laws are those
probability measures. These asymptotic bounds are expressed in terms of a \emph{speed function} $v_t$ (that tends to infinity)
and a lower semicontinuous \emph{rate function} $I:\mathcal{X}\to[0,\infty]$. The concept of large deviation principle is a 
basic definition in the theory of large deviations; this theory allows us to compute the probabilities of rare events on an 
exponential scale (see \cite{DemboZeitouni} as a reference of this topic).

The term \emph{moderate deviations} is often used in the literature to mean a class of large deviation principles 
that, in some sense, fills the gap between two asymptotic regimes:
\begin{enumerate}
\item the convergence of $C_t$ in probability to zero, which is governed by a large deviation principle with speed $v_t$; 
\item the weak convergence of $\sqrt{v_t} C_t$ to a centered Normal distribution.
\end{enumerate}
The speed functions and the random variables involved in these large deviation principles depend on some scalings in a suitable
class; moreover, the large deviation principles in this class are governed by the same quadratic rate function that uniquely 
vanishes at zero. Typically the scalings consist of families of positive numbers $\{a_t:t>0\}$ such that
$$a_t\to 0\ \mbox{and}\ v_ta_t\to\infty,$$ 
and one can show that $\{\sqrt{v_t a_t}C_t\}$ satisfies the large deviation principle with speed $\frac{1}{a_t}$; note that 
$\frac{1}{a_t}$ has a lower intensity than the speed $v_t$, and this explains the use of the term \emph{moderate}. We also recall
that we recover the two asymptotic regimes stated above for $a_t=\frac{1}{v_t}$ (in this case $v_ta_t\to\infty$ fails) and for 
$a_t=1$ (in this case $a_t\to 0$ fails).

The term \emph{noncentral moderate deviations} has been recently used in the literature when we have a class of large deviation 
principles that, in some sense, fills the gap between a convergence to a constant (typically zero) and the weak convergence
towards a non-Gaussian distribution. Some examples of noncentral moderate deviations can be found in \cite{GiulianoMacci}, where 
the weak convergences are towards Gumbel, exponential, and Laplace distributions. In that reference, the interested reader can find 
some other previous references in the literature with some other examples.

The aim of this paper is to present some examples of noncentral moderate deviations based on \emph{fractional Skellam
processes}. In these examples we always have $v_t=t$, and therefore the scalings are families of positive numbers 
$\{a_t:t>0\}$ such that
\begin{equation}\label{eq:MDconditions}
a_t\to 0\ \mbox{and}\ ta_t\to\infty.
\end{equation}
Skellam processes are given by the difference of two independent Poisson processes. In this paper we consider two
fractional Skellam processes studied in \cite{KerssLeonenkoSikorskii}; some more recent generalized versions of these 
processes can be found in \cite{GuptaKumarLeonenko} and \cite{KatariaKhandakar-arxiv}. The fractional Skellam processes in 
\cite{KerssLeonenkoSikorskii} are closely related to the definition of the fractional Poisson process in the literature.
We recall that a fractional Poisson process is obtained as an independent random time-change of a Poisson process with 
an inverse of stable subordinator (see e.g. \cite{BeghinOrsingher2009}, \cite{BeghinOrsingher2010} and 
\cite{MeerschaertNaneVellaisamy}); here we are referring to the time fractional Poisson process and, for the definitions 
of space and space-time fractional Poisson process, the interested reader can refer to \cite{OrsingherPolito} (see also
\cite{KatariaKhandakarSAA2022} as a very recent paper on time-changed space-time fractional Poisson processes). Then the 
fractional Skellam processes studied in \cite{KerssLeonenkoSikorskii} are obtained in a quite natural way as follows.
\begin{itemize}
\item \emph{Fractional Skellam process of type 1}. A difference between two independent fractional Poisson processes 
(so we have two independent random time-changes for each one of the involved fractional Poisson processes);
\item \emph{Fractional Skellam process of type 2}. An independent random time-change of a Skellam process with an inverse 
of stable subordinator.
\end{itemize}
It is easy to check (see Remark \ref{rem:CFPP-particular-case}) that the fractional Skellam process of type 2 is a particular 
compound fractional Poisson process (and this is not surprising because a Skellam process is a particular compound Poisson 
process; see Remark \ref{rem:CPP-particular-case}). Therefore, the moderate deviation results for the fractional Skellam process
of type 2 can be obtained from the ones in \cite{BeghinMacciSPL2022}. Here, since we deal with random 
time-changes of Skellam processes, for completeness we recall the references \cite{BuchakSakhno2017}, 
\cite{BuchakSakhno2019} and \cite{KatariaKhandakarFCAA2022}.

Here for completeness we present a brief review of the references with results on large/moderate deviations
for fractional Poisson processes or similar models: \cite{BeghinMacciSPL2013} and \cite{BeghinMacciSPL2017} with results for the
(possibly multivariate) alternative fractional Poisson process, \cite{BeghinMacciMartinucci} with results for random time-changed
continuous-time Markov chains on integers with alternating rates, \cite{MacciPacchiarottiVilla} with results for a non-standard
model based on the Prabhakar function in \cite{PoganyTomovski}, and for a state dependent model in \cite{GarraOrsingherPolito}.

We conclude with the outline of the paper. We start with some preliminaries in Section \ref{sec:preliminaries}.
The results for the fractional Skellam processes of type 1 and 2 are presented in Sections \ref{sec:type1} and
\ref{sec:type2}, respectively. In Section \ref{sec:comparisons} we compare some rate functions and  
present some plots. Finally, in Section \ref{sec:concluding-remarks}, we present some concluding remarks.

\section{Preliminaries}\label{sec:preliminaries}
In this section, we recall some preliminaries on large deviations and on fractional Skellam processes.

\subsection{On large deviations}
Here we present definitions and results for families of real random variables $\{Z_t:t>0\}$ defined on 
the same probability space $(\Omega,\mathcal{F},P)$; moreover, in view of what follows, we consider the case $t\to\infty$.
We start with the definition of large deviation principle (see e.g. \cite{DemboZeitouni}, pages 4-5). A family of numbers 
$\{v_t:t>0\}$ such that $v_t\to\infty$ (as $t\to\infty$) is called a \emph{speed function}, and a lower semicontinuous 
function $I:\mathbb{R}\to[0,\infty]$ is called a \emph{rate function}. Then $\{Z_t:t>0\}$ satisfies the large deviation 
principle (LDP from now on) with speed $v_t$ and a rate function $I$ if
$$\limsup_{t\to\infty}\frac{1}{v_t}\log P(Z_t\in C)\leq-\inf_{x\in C}I(x)\quad\mbox{for all closed sets}\ C,$$
and
$$\liminf_{t\to\infty}\frac{1}{v_t}\log P(Z_t\in O)\geq-\inf_{x\in O}I(x)\quad\mbox{for all open sets}\ O.$$
Moreover the rate function $I$ is said to be \emph{good} if, for every $\beta\geq 0$, the level set 
$\{x\in\mathbb{R}:I(x)\leq\beta\}$ is compact. The following well-known theorem provides a sufficient condition to have 
an LDP, and makes it easy to compute the corresponding speed and rate functions (see e.g. Theorem 2.3.6(c) in 
\cite{DemboZeitouni}).

\begin{theorem}[G\"artner Ellis Theorem]\label{th:GE}
Assume that, for all $\theta\in\mathbb{R}$, there exists
$$\Lambda(\theta):=\lim_{t\to\infty}\frac{1}{v_t}\log\mathbb{E}\left[e^{v_t\theta Z_t}\right]$$
as an extended real number; moreover assume that the origin $\theta=0$ 
belongs to the interior of the set
$$\mathcal{D}(\Lambda):=\{\theta\in\mathbb{R}:\Lambda(\theta)<\infty\}.$$
Furthermore let $\Lambda^*$ be the Legendre-Fenchel transform of the function $\Lambda$, i.e.
the function defined by
$$\Lambda^*(x):=\sup_{\theta\in\mathbb{R}}\{\theta x-\Lambda(\theta)\}.$$
Then, if $\Lambda$ is essentially smooth and lower semi-continuous, then $\{Z_t:t>0\}$
satisfies the LDP with speed $v_t$ and good rate function $\Lambda^*$.
\end{theorem}

We also recall (see e.g. Definition 2.3.5 in \cite{DemboZeitouni}) that $\Lambda$ is essentially smooth
if the interior of $\mathcal{D}(\Lambda)$ is non-empty, the function $\Lambda$ is differentiable 
throughout the interior of $\mathcal{D}(\Lambda)$, and $\Lambda$ is steep, i.e. $|\Lambda^\prime(\theta_n)|\to\infty$
whenever $\{\theta_n:n\geq 1\}$ is a sequence of points in the interior of $\mathcal{D}(\Lambda)$ which converge to 
a boundary point of $\mathcal{D}(\Lambda)$. A particular simple case (which always occurs in the applications
of the G\"artner Ellis Theorem in this paper) is when $\mathcal{D}(\Lambda)=\mathbb{R}$ and $\Lambda$ is a
differentiable function; indeed, in such a case, the function $\Lambda$ is essentially smooth (the steepness 
condition holds vacuously) and lower semi-continuous.

\subsection{On fractional Skellam processes}
We start with the definition of the (non-fractional) Skellam process. Let $\{N_{\lambda_1}(t):t\geq 0\}$ 
and $\{N_{\lambda_2}(t):t\geq 0\}$ be two independent Poisson processes with intensities $\lambda_1>0$ 
and $\lambda_2>0$, respectively. In particular we consider the notation $\underline{\lambda}=(\lambda_1,\lambda_2)$.
Then the process $\{S_{\underline{\lambda}}(t):t\geq 0\}$ defined by
$$S_{\underline{\lambda}}(t):=N_{\lambda_1}(t)-N_{\lambda_2}(t)$$
is called Skellam process. Moreover, for each fixed $t\geq 0$, we have
$$\mathbb{E}[e^{\theta S_{\underline{\lambda}}(t)}]
=\exp((\lambda_1(e^\theta-1)+\lambda_2(e^{-\theta}-1))t)\ (\mbox{for all}\ \theta\in\mathbb{R}).$$

\begin{remark}\label{rem:CPP-particular-case}
It is easy to check that $\{S_{\underline{\lambda}}(t):t\geq 0\}$ can be seen
as a compound Poisson process $\{\sum_{k=1}^{N_{\lambda_1+\lambda_2}(t)}X_k:t\geq 0\}$,
where $\{X_k:k\geq 1\}$ and $\{N_{\lambda_1+\lambda_2}(t):t\geq 0\}$ are independent,
$\{X_k:k\geq 1\}$ are i.i.d. random variables such that
$$P(X_k=1)=1-P(X_k=-1)=\frac{\lambda_1}{\lambda_1+\lambda_2},$$
and $\{N_{\lambda_1+\lambda_2}(t):t\geq 0\}$ is a Poisson process with intensity
$\lambda_1+\lambda_2$.
\end{remark}

In view of what follows we recall some other preliminaries. We start with the definition of the Mittag-Leffler function
(see e.g. \cite{GorenfloKilbasMainardiRogosin}, eq. (3.1.1))
$$E_\nu(x):=\sum_{k=0}^\infty\frac{x^k}{\Gamma(\nu k+1)}\quad\mbox{for}\ \nu,x\in\mathbb{C}.$$
Actually throughout this paper we have $\nu\in(0,1)$ and $x\in\mathbb{R}$; moreover it is known (see 
Proposition 3.6 in \cite{GorenfloKilbasMainardiRogosin} for the case $\alpha\in(0,2)$; indeed $\alpha$ in that reference
coincides with $\nu$ in this paper) that we have
\begin{equation}\label{eq:ML-asymptotics}
\left\{\begin{array}{l}
E_\nu(x)\sim\frac{e^{x^{1/\nu}}}{\nu},\ \mbox{as}\ x\to\infty\\
\frac{1}{x}\log E_\nu(x)\to 0\ \mbox{as}\ x\to-\infty
\end{array}\right.
\end{equation}
(this is the correct version of eq. (3) in \cite{BeghinMacciSPL2022}; indeed we need the condition presented here for 
$x\to-\infty$, instead of $E_\nu(x)\to 0$).

Now we recall some moment generating functions which can be expressed in terms of the Mittag-Leffler
function. If we consider the inverse of the stable subordinator $\{L_\nu(t):t\geq 0\}$, then we have 
\begin{equation}\label{eq:MGF-inverse-stable-sub}
\mathbb{E}[e^{\theta L_\nu(t)}]=E_\nu(\theta t^\nu)\ (\mbox{for all}\ \theta\in\mathbb{R}).
\end{equation}
This formula appears in several references with $\theta\leq 0$ only; however, this restriction is not
needed because we can refer to the analytic continuation of the Laplace transform with complex argument.

The fractional Poisson process $\{N_{\nu,\lambda}(t):t\geq 0\}$ is defined by
$$N_{\nu,\lambda}(t):=N_\lambda(L_\nu(t)),$$
where $\{N_\lambda(t):t\geq 0\}$ is a (non-fractional) Poisson process with intensity $\lambda$,
independent of $\{L_\nu(t):t\geq 0\}$; it is known that
$$\mathbb{E}[e^{\theta N_{\nu,\lambda}(t)}]=E_\nu(\lambda(e^\theta-1)t^\nu)\ (\mbox{for all}\ \theta\in\mathbb{R}).$$

Now we are ready to provide the definitions of two fractional Skellam processes and their moment generating functions 
(see \cite{KerssLeonenkoSikorskii}, Definitions 3.1-3.2 and Theorems 3.1-3.2). In particular we consider the notation 
$\underline{\nu}=(\nu_1,\nu_2)$ for $\nu_1,\nu_2\in(0,1)$.

\paragraph{Fractional Skellam process of type 1.}
It is the process $\{Y_{\underline{\nu},\underline{\lambda}}(t):t\geq 0\}$ defined by
$$Y_{\underline{\nu},\underline{\lambda}}(t):=N_{\nu_1,\lambda_1}(t)-N_{\nu_2,\lambda_2}(t),$$
where $\{N_{\nu_1,\lambda_1}(t):t\geq 0\}$ and $\{N_{\nu_2,\lambda_2}(t):t\geq 0\}$ are two independent
fractional Poisson processes. Then we have
$$\mathbb{E}[e^{\theta Y_{\underline{\nu},\underline{\lambda}}(t)}]
=E_{\nu_1}(\lambda_1(e^\theta-1)t^{\nu_1})E_{\nu_2}(\lambda_2(e^{-\theta}-1)t^{\nu_2})
\ (\mbox{for all}\ \theta\in\mathbb{R}).$$

\begin{remark}\label{rem:simplified-notation}
Some results for $\{Y_{\underline{\nu},\underline{\lambda}}(t):t\geq 0\}$ presented below concern
the case $\nu_1=\nu_2$; in this case we set $\nu_1=\nu_2=\nu$ for some $\nu\in(0,1)$, and use the
notation $Y_{\nu,\underline{\lambda}}(t)$ in place of $Y_{\underline{\nu},\underline{\lambda}}(t)$.
\end{remark}

\paragraph{Fractional Skellam process of type 2.}
It is the process $\{Z_{\nu,\underline{\lambda}}(t):t\geq 0\}$ defined by
$$Z_{\nu,\underline{\lambda}}(t):=S_{\underline{\lambda}}(L_\nu(t)),$$
where the Skellam process $\{S_{\underline{\lambda}}(t):t\geq 0\}$ and the inverse of the
stable subordinator $\{L_\nu(t):t\geq 0\}$ are independent. Then we have
$$\mathbb{E}[e^{\theta Z_{\nu,\underline{\lambda}}(t)}]
=E_\nu((\lambda_1(e^\theta-1)+\lambda_2(e^{-\theta}-1))t^\nu)\ (\mbox{for all}\ \theta\in\mathbb{R}).$$

\begin{remark}\label{rem:CFPP-particular-case}
We have recalled that the fractional Poisson process can be seen as a time changed (non-fractional)
Poisson process with an independent inverse of the stable subordinator. Then, by taking into account
Remark \ref{rem:CPP-particular-case}, we can say that $\{Z_{\nu,\underline{\lambda}}(t):t\geq 0\}$ 
is distributed as the compound fractional Poisson process $\{\sum_{k=1}^{N_{\nu,\lambda_1+\lambda_2}(t)}X_k:t\geq 0\}$,
where $\{X_k:k\geq 1\}$ and $\{N_{\nu,\lambda_1+\lambda_2}(t):t\geq 0\}$ are independent, $\{X_k:k\geq 1\}$
are i.i.d. random variables as in Remark \ref{rem:CPP-particular-case}, and $\{N_{\nu,\lambda_1+\lambda_2}(t):t\geq 0\}$ 
is a fractional Poisson process.
\end{remark}

\begin{remark}\label{rem:symmetry}
Assume that $\nu_1=\nu_2=\nu$ for some $\nu\in(0,1)$ (and recall the slight change of notation explained in Remark 
\ref{rem:simplified-notation} for the fractional Skellam process of type 1). Then, if $\lambda_1=\lambda_2$, the random 
variables $Y_{\nu,\underline{\lambda}}(t)$ and $Z_{\nu,\underline{\lambda}}(t)$ are symmetric (around zero); namely 
$Y_{\nu,\underline{\lambda}}(t)$ and $Z_{\nu,\underline{\lambda}}(t)$ are distributed as $-Y_{\nu,\underline{\lambda}}(t)$ 
and $-Z_{\nu,\underline{\lambda}}(t)$, respectively. Then we have some consequences highlighted in Remarks 
\ref{rem:symmetry-type1} and \ref{rem:symmetry-type2}.
\end{remark}

\section{Noncentral moderate deviations for the type 1 process}\label{sec:type1}
We start with the first result for which we could have $\nu_1\neq \nu_2$.

\begin{proposition}\label{prop:LD-type1}
Let $\Psi_{\underline{\nu},\underline{\lambda}}^{(1)}$ be the function defined by
\begin{equation}\label{eq:LD-GE-limit-type1}
\Psi_{\underline{\nu},\underline{\lambda}}^{(1)}(\theta):=\left\{\begin{array}{ll}
(\lambda_1(e^\theta-1))^{1/\nu_1}&\ \mbox{if}\ \theta\geq 0\\
(\lambda_2(e^{-\theta}-1))^{1/\nu_2}&\ \mbox{if}\ \theta<0.
\end{array}\right.
\end{equation}
Then $\left\{\frac{Y_{\underline{\nu},\underline{\lambda}}(t)}{t}:t>0\right\}$ satisfies the LDP with 
speed $v_t=t$ and good rate function $I_{\mathrm{LD}}^{(1)}$ defined by
\begin{equation}\label{eq:LD-type1-rf-Legendre-transform}
I_{\mathrm{LD}}^{(1)}(x):=\sup_{\theta\in\mathbb{R}}\{\theta x-\Psi_{\underline{\nu},\underline{\lambda}}^{(1)}(\theta)\}.
\end{equation}

\end{proposition}
\begin{proof}
We prove this proposition by applying the G\"artner Ellis Theorem. More precisely we have to show that
\begin{equation}\label{eq:LD-GET-limit-type1}
\lim_{t\to\infty}\frac{1}{t}\log\mathbb{E}\left[e^{t\theta \frac{Y_{\underline{\nu},\underline{\lambda}}(t)}{t}}\right]
=\Psi_{\underline{\nu},\underline{\lambda}}^{(1)}(\theta)\ (\mbox{for all}\ \theta\in\mathbb{R}),
\end{equation}
where $\Psi_{\underline{\nu},\underline{\lambda}}^{(1)}$ is the function in \eqref{eq:LD-GE-limit-type1}.

The case $\theta=0$ is immediate. Moreover, we remark that
$$\log\mathbb{E}\left[e^{t\theta \frac{Y_{\underline{\nu},\underline{\lambda}}(t)}{t}}\right]
=\log\mathbb{E}\left[e^{\theta Y_{\underline{\nu},\underline{\lambda}}(t)}\right]
=\log E_{\nu_1}(\lambda_1(e^\theta-1)t^{\nu_1})+\log E_{\nu_2}(\lambda_2(e^{-\theta}-1)t^{\nu_2}).$$
Then, by taking into account the asymptotic behaviour of the Mittag-Leffler function in \eqref{eq:ML-asymptotics}, we have
$$\lim_{t\to\infty}\frac{1}{t}\log E_{\nu_1}(\lambda_1(e^\theta-1)t^{\nu_1})+\lim_{t\to\infty}
\frac{1}{t}\log E_{\nu_2}(\lambda_2(e^{-\theta}-1)t^{\nu_2})=(\lambda_1(e^\theta-1))^{1/\nu_1}\ \mbox{for}\ \theta>0,$$
and
$$\lim_{t\to\infty}\frac{1}{t}\log E_{\nu_1}(\lambda_1(e^\theta-1)t^{\nu_1})+\lim_{t\to\infty}
\frac{1}{t}\log E_{\nu_2}(\lambda_2(e^{-\theta}-1)t^{\nu_2})=(\lambda_2(e^{-\theta}-1))^{1/\nu_2}\ \mbox{for}\ \theta<0.$$
Thus the limit in \eqref{eq:LD-GET-limit-type1} is checked.

In conclusion the desired LDP holds noting that the function $\Psi_{\underline{\nu},\underline{\lambda}}^{(1)}$
in \eqref{eq:LD-GE-limit-type1} is finite (for all $\theta\in\mathbb{R}$) and differentiable.
\end{proof}

\begin{remark}\label{rem:LD-type1-nu1=nu2=1/2}
In general, we do not have an explicit expression for the rate function $I_{\mathrm{LD}}^{(1)}$ in Proposition \ref{prop:LD-type1}
(see \eqref{eq:LD-type1-rf-Legendre-transform}). However, if $\nu_1=\nu_2=\frac{1}{2}$, then we have
$$I_{\mathrm{LD}}^{(1)}(x)=\left\{\begin{array}{ll}
x\log\left(\frac{1}{2}+\frac{1}{2}\sqrt{1+\frac{2x}{\lambda_1^2}}\right)-\left(\frac{1}{2}\sqrt{\lambda_1^2+2x}-\frac{\lambda_1}{2}\right)^2&\ \mbox{if}\ x\geq 0\\
-x\log\left(\frac{1}{2}+\frac{1}{2}\sqrt{1-\frac{2x}{\lambda_2^2}}\right)-\left(\frac{1}{2}\sqrt{\lambda_2^2-2x}-\frac{\lambda_2}{2}\right)^2&\ \mbox{if}\ x<0.
\end{array}\right.$$
Indeed, after some computations, one can check that the supremum in \eqref{eq:LD-type1-rf-Legendre-transform} 
(with $\nu_1=\nu_2=\frac{1}{2}$) is attained at $\theta=0$ for $x=0$, at 
$\theta=\log\left(\frac{1}{2}+\frac{1}{2}\sqrt{1+\frac{2x}{\lambda_1^2}}\right)$ for $x>0$,
and at $\theta=-\log\left(\frac{1}{2}+\frac{1}{2}\sqrt{1-\frac{2x}{\lambda_2^2}}\right)$ for $x<0$.
\end{remark}

From now on we assume that $\nu_1$ and $\nu_2$ coincide, and therefore we consider the change of notation
in Remark \ref{rem:simplified-notation} for $\nu=\nu_1=\nu_2$. Moreover we set
\begin{equation}\label{eq:exponent}
\alpha_1(\nu):=1-\nu.
\end{equation}

\begin{proposition}\label{prop:weak-convergence-type1}
Assume that $\nu_1=\nu_2=\nu$ for some $\nu\in(0,1)$ and let $\alpha_1(\nu)$ be defined in \eqref{eq:exponent}. 
Then $\{t^{\alpha_1(\nu)}\frac{Y_{\nu,\underline{\lambda}}(t)}{t}:t>0\}$ converges weakly to 
$\lambda_1L_\nu^\circ(1)-\lambda_2L_\nu^{\circ\circ}(1)$, where $L_\nu^\circ(1)$ and $L_\nu^{\circ\circ}(1)$ are two i.i.d. 
random variables distributed as $L_\nu(1)$.
\end{proposition}
\begin{proof}
We have to check that
$$\lim_{t\to\infty}\mathbb{E}\left[e^{\theta t^{\alpha_1(\nu)}\frac{Y_{\nu,\underline{\lambda}}(t)}{t}}\right]=
\underbrace{\mathbb{E}[e^{\theta (\lambda_1L_\nu^\circ(1)-\lambda_2L_\nu^{\circ\circ}(1))}]}
_{=E_\nu(\lambda_1\theta)E_\nu(-\lambda_2\theta)}\ (\mbox{for all}\ \theta\in\mathbb{R})$$
(here we take into account that $L_\nu^\circ(1)$ and $L_\nu^{\circ\circ}(1)$ are i.i.d., and the expression of the
moment generating function in eq. \eqref{eq:MGF-inverse-stable-sub}). This can be readily done noting that, 
for two suitable remainders $o\left(\frac{1}{t^{2\nu}}\right)$ such that $t^{2\nu}o\left(\frac{1}{t^{2\nu}}\right)\to 0$,
we have
\begin{multline*}
\mathbb{E}\left[e^{\theta t^{\alpha_1(\nu)}\frac{Y_{\nu,\underline{\lambda}}(t)}{t}}\right]
=\mathbb{E}\left[e^{\theta\frac{Y_{\nu,\underline{\lambda}}(t)}{t^\nu}}\right]
=E_\nu(\lambda_1(e^{\theta/t^\nu}-1)t^\nu)E_\nu(\lambda_2(e^{-\theta/t^\nu}-1)t^\nu)\\
=E_\nu\left(\lambda_1\left(\frac{\theta}{t^\nu}+\frac{\theta^2}{2t^{2\nu}}+
o\left(\frac{1}{t^{2\nu}}\right)\right)t^\nu\right)E_\nu\left(\lambda_2\left(-\frac{\theta}{t^\nu}+\frac{\theta^2}{2t^{2\nu}}+
o\left(\frac{1}{t^{2\nu}}\right)\right)t^\nu\right)\\
=E_\nu\left(\lambda_1\left(\theta+\frac{\theta^2}{2t^\nu}+
t^\nu o\left(\frac{1}{t^{2\nu}}\right)\right)\right)E_\nu\left(\lambda_2\left(-\theta+\frac{\theta^2}{2t\nu}+
t^\nu o\left(\frac{1}{t^{2\nu}}\right)\right)\right);
\end{multline*}
then we get the desired limit letting $t$ go to infinity (for each fixed $\theta\in\mathbb{R}$).
\end{proof}

\begin{proposition}\label{prop:NCMD-type1}
Assume that $\nu_1=\nu_2=\nu$ for some $\nu\in(0,1)$ and let $\alpha_1(\nu)$ be defined in \eqref{eq:exponent}. 
Then, for every family of positive numbers $\{a_t:t>0\}$ such that \eqref{eq:MDconditions} holds, the family 
of random variables $\left\{\frac{(a_tt)^{\alpha_1(\nu)}Y_{\nu,\underline{\lambda}}(t)}{t}:t>0\right\}$ satisfies
the LDP with speed $1/a_t$ and good rate function $I_{\mathrm{MD}}^{(1)}$ defined by
$$I_{\mathrm{MD}}^{(1)}(x):=\left\{\begin{array}{ll}
(\nu^{\nu/(1-\nu)}-\nu^{1/(1-\nu)})\left(\frac{x}{\lambda_1}\right)^{1/(1-\nu)}&\ \mbox{if}\ x\geq 0\\
(\nu^{\nu/(1-\nu)}-\nu^{1/(1-\nu)})\left(\frac{x}{-\lambda_2}\right)^{1/(1-\nu)}&\ \mbox{if}\ x<0.
\end{array}\right.$$
\end{proposition}
\begin{proof}
We prove this proposition by applying the G\"artner Ellis Theorem. More precisely we have to show that
\begin{equation}\label{eq:NCMD-GET-limit-type1}
\lim_{t\to\infty}\frac{1}{1/a_t}\log\mathbb{E}\left[e^{\frac{\theta}{a_t}\frac{(a_tt)^{\alpha_1(\nu)}Y_{\nu,\underline{\lambda}}(t)}{t}}\right]
=\tilde{\Psi}_{\nu,\underline{\lambda}}^{(1)}(\theta)\ (\mbox{for all}\ \theta\in\mathbb{R}),
\end{equation}
where $\tilde{\Psi}_{\nu,\underline{\lambda}}^{(1)}$ is the function defined by
$$\tilde{\Psi}_{\nu,\underline{\lambda}}^{(1)}(\theta):=\left\{\begin{array}{ll}
(\lambda_1\theta)^{1/\nu}&\ \mbox{if}\ \theta\geq 0\\
(-\lambda_2\theta)^{1/\nu}&\ \mbox{if}\ \theta<0;
\end{array}\right.$$
indeed, since the function $\tilde{\Psi}_{\nu,\underline{\lambda}}^{(1)}$ is finite (for 
all $\theta\in\mathbb{R}$) and differentiable, the desired LDP holds noting that the Legendre-Fenchel
transform $(\tilde{\Psi}_{\nu,\underline{\lambda}}^{(1)})^*$ of $\tilde{\Psi}_{\nu,\underline{\lambda}}^{(1)}$, i.e.
the function $(\tilde{\Psi}_{\nu,\underline{\lambda}}^{(1)})^*$ defined by
$$(\tilde{\Psi}_{\nu,\underline{\lambda}}^{(1)})^*(x):=
\sup_{\theta\in\mathbb{R}}\{\theta x-\tilde{\Psi}_{\nu,\underline{\lambda}}^{(1)}(\theta)\}
\ (\mbox{for all}\ x\in\mathbb{R}),$$
coincides with the function $I_{\mathrm{MD}}^{(1)}$ in the statement of the proposition (for $x=0$ the
supremum is attained at $\theta=0$, for $x>0$ the supremum is attained at 
$\theta=\frac{1}{\lambda_1}(\frac{\nu x}{\lambda_1})^{\nu/(1-\nu)}$, for $x<0$ the supremum is
attained at $\theta=-\frac{1}{\lambda_2}(\frac{\nu x}{-\lambda_2})^{\nu/(1-\nu)}$).

So we conclude the proof by checking the limit in \eqref{eq:NCMD-GET-limit-type1}. The case $\theta=0$ 
is immediate. Moreover we remark that, for two suitable remainders $o\left(\frac{1}{(a_tt)^{2\nu}}\right)$
such that $(a_tt)^{2\nu}o\left(\frac{1}{(a_tt)^{2\nu}}\right)\to 0$, we have
\begin{multline*}
\log\mathbb{E}\left[e^{\frac{\theta}{a_t}\frac{(a_tt)^{\alpha_1(\nu)}Y_{\nu,\underline{\lambda}}(t)}{t}}\right]
=\log\mathbb{E}\left[e^{\theta\frac{Y_{\nu,\underline{\lambda}}(t)}{(a_tt)^\nu}}\right]\\
=\log E_\nu(\lambda_1(e^{\theta/(a_tt)^\nu}-1)t^\nu)+\log E_\nu(\lambda_2(e^{-\theta/(a_tt)^\nu}-1)t^\nu)\\
=\log E_\nu\left(\lambda_1\left(\frac{\theta}{(a_tt)^\nu}+\frac{\theta^2}{2(a_tt)^{2\nu}}+
o\left(\frac{1}{(a_tt)^{2\nu}}\right)\right)t^\nu\right)\\
+\log E_\nu\left(\lambda_2\left(-\frac{\theta}{(a_tt)^\nu}+\frac{\theta^2}{2(a_tt)^{2\nu}}+
o\left(\frac{1}{(a_tt)^{2\nu}}\right)\right)t^\nu\right)\\
=\log E_\nu\left(\frac{\lambda_1}{a_t^\nu}\left(\theta+\frac{\theta^2}{2(a_tt)^\nu}+
(a_tt)^\nu o\left(\frac{1}{(a_tt)^{2\nu}}\right)\right)\right)\\
+\log E_\nu\left(\frac{\lambda_2}{a_t^\nu}\left(-\theta+\frac{\theta^2}{2(a_tt)^\nu}+
(a_tt)^\nu o\left(\frac{1}{(a_tt)^{2\nu}}\right)\right)\right).
\end{multline*}
Then, by taking into account the asymptotic behaviour of the Mittag-Leffler function in \eqref{eq:ML-asymptotics},
we have
$$\lim_{t\to\infty}\frac{1}{1/a_t}\log \mathbb{E}\left[e^{\frac{\theta}{a_t}\frac{(a_tt)^{\alpha_1(\nu)}Y_{\nu,\underline{\lambda}}(t)}{t}}\right]=(\lambda_1\theta)^{1/\nu}\ \mbox{for}\ \theta>0,$$
and
$$\lim_{t\to\infty}\frac{1}{1/a_t}\log \mathbb{E}\left[e^{\frac{\theta}{a_t}\frac{(a_tt)^{\alpha_1(\nu)}Y_{\nu,\underline{\lambda}}(t)}{t}}\right]=(-\lambda_2\theta)^{1/\nu}\ \mbox{for}\ \theta<0.$$
Thus the limit in \eqref{eq:NCMD-GET-limit-type1} is checked.
\end{proof}

\begin{remark}\label{rem:weakconv-NCMD-type1}
The set $\{x\in\mathbb{R}:I_{\mathrm{MD}}^{(1)}(x)<\infty\}=\mathbb{R}$ (see Proposition \ref{prop:NCMD-type1}) coincides with 
the support of the weak limit in Proposition \ref{prop:weak-convergence-type1}.
\end{remark}

\begin{remark}\label{rem:symmetry-type1}
Assume that $\lambda_1=\lambda_2$. Then: if $\nu_1=\nu_2$, the rate function $I_{\mathrm{LD}}^{(1)}$ in Proposition \ref{prop:LD-type1}
is a symmetric function (we can say this, even if an explicit expression of $I_{\mathrm{LD}}^{(1)}$ is not available, because $\Psi_{\underline{\nu},\underline{\lambda}}^{(1)}$ is a symmetric function); the weak limit in Proposition \ref{prop:weak-convergence-type1}
is a symmetric random variable; the rate function $I_{\mathrm{MD}}^{(1)}$ in Proposition \ref{prop:NCMD-type1} is a symmetric function.
\end{remark}

\section{Noncentral moderate deviations for the type 2 process}\label{sec:type2}
The results in this section can be derived directly from the results in \cite{BeghinMacciSPL2022}; indeed, by taking into account Remark 
\ref{rem:CFPP-particular-case}, the fractional Skellam process of type 2 is a particular compound fractional Poisson process. So we only
give the statements of propositions without proofs.

\begin{proposition}\label{prop:LD-type2}
Let $\Psi_{\nu,\underline{\lambda}}^{(2)}$ be the function defined by
\begin{equation}\label{eq:LD-GE-limit-type2}
\Psi_{\nu,\underline{\lambda}}^{(2)}(\theta):=\left\{\begin{array}{ll}
(\lambda_1(e^\theta-1)+\lambda_2(e^{-\theta}-1))^{1/\nu}&\ \mbox{if}\ \lambda_1(e^\theta-1)+\lambda_2(e^{-\theta}-1)\geq 0\\
0&\ \mbox{if}\ \lambda_1(e^\theta-1)+\lambda_2(e^{-\theta}-1)<0.
\end{array}\right.
\end{equation}
Then $\left\{\frac{Z_{\nu,\underline{\lambda}}(t)}{t}:t>0\right\}$ satisfies the LDP with 
speed $v_t=t$ and good rate function $I_{\mathrm{LD}}^{(2)}$ defined by
$$I_{\mathrm{LD}}^{(2)}(x):=\sup_{\theta\in\mathbb{R}}\{\theta x-\Psi_{\nu,\underline{\lambda}}^{(2)}(\theta)\}.$$
\end{proposition}

\begin{remark}\label{rem:MSTA2021}
The LDP in Proposition \ref{prop:LD-type2} with $\lambda_1=\lambda_2=\lambda$ coincides with the LDP in 
Proposition 4.2 in \cite{BeghinMacciMartinucci} with $\alpha_1=\alpha_2=\beta_1=\beta_2=\lambda$. Indeed the 
function $\Lambda_\nu$ in that proposition is defined by
$$\Lambda_\nu(\theta):=\left\{\begin{array}{ll}
(\Lambda(\theta))^{1/\nu}&\ \mbox{if}\ \Lambda(\theta)\geq 0\\
0&\ \mbox{if}\ \Lambda(\theta)<0,
\end{array}\right.$$
where, after some computations, one can check that
$$\Lambda(\theta)=\frac{\lambda(e^{2\theta}+1)}{e^\theta}-2\lambda=\lambda(e^\theta-1+e^{-\theta}-1).$$
\end{remark}
The next two propositions can be derived from Propositions 3.2-3.3 in \cite{BeghinMacciSPL2022}. By Remark 
\ref{rem:CFPP-particular-case} we have $Z_{\nu,\underline{\lambda}}(t)\stackrel{d}{=}\sum_{k=1}^{N_{\nu,\lambda_1+\lambda_2}(t)}X_k$,
where
$$\mu:=E(X_k)=\frac{\lambda_1-\lambda_2}{\lambda_1+\lambda_2}\ \mbox{and}\ 
\sigma^2:=\mathrm{Var}(X_k)=\frac{4\lambda_1\lambda_2}{(\lambda_1+\lambda_2)^2}$$
coincide with $\mu$ and $\sigma^2$ in \cite{BeghinMacciSPL2022}. Moreover we set
\begin{equation}\label{eq:exponents}
\alpha_2(\nu):=\left\{\begin{array}{ll}
1-\nu/2&\ \mbox{if}\ \mu=0,\ \mbox{i.e. if}\ \lambda_1=\lambda_2,\\
1-\nu&\ \mbox{if}\ \mu\neq 0,\ \mbox{i.e. if}\ \lambda_1\neq\lambda_2,
\end{array}\right.
\end{equation}
which coincides with $\alpha(\nu)$ in \cite{BeghinMacciSPL2022}. Also note that, if $\lambda_1=\lambda_2=\lambda$ for some $\lambda>0$,
we take into account that $\sqrt{\frac{4\lambda_1\lambda_2}{\lambda_1+\lambda_2}}=\sqrt{2\lambda}$ in Proposition 
\ref{prop:weak-convergence-type2}, and $\frac{\lambda_1+\lambda_2}{2\lambda_1\lambda_2}=\frac{1}{\lambda}$ in Proposition 
\ref{prop:NCMD-type2}.

\begin{proposition}\label{prop:weak-convergence-type2}
Let $\alpha_2(\nu)$ be defined in \eqref{eq:exponents}. Then:
\begin{itemize}
\item if $\lambda_1=\lambda_2=\lambda$ for some $\lambda>0$, then $\{t^{\alpha_2(\nu)}\frac{Z_{\nu,\underline{\lambda}}(t)}{t}:t>0\}$
converges weakly to $\sqrt{2\lambda L_\nu(1)}W$, where $W$ is a standard Normal distributed random variable, and independent to 
$L_\nu(1)$;
\item if $\lambda_1\neq\lambda_2$, then $\{t^{\alpha_2(\nu)}\frac{Z_{\nu,\underline{\lambda}}(t)}{t}:t>0\}$ converges weakly to
$(\lambda_1-\lambda_2)L_\nu(1)$.
\end{itemize}
\end{proposition}

\begin{proposition}\label{prop:NCMD-type2}
Let $\alpha_2(\nu)$ be defined in \eqref{eq:exponents}. Then, for every family of positive numbers $\{a_t:t>0\}$ 
such that \eqref{eq:MDconditions} holds, the family of random variables 
$\left\{\frac{(a_tt)^{\alpha_2(\nu)}Z_{\nu,\underline{\lambda}}(t)}{t}:t>0\right\}$ satisfies the LDP with speed $1/a_t$ and 
good rate function $I_{\mathrm{MD},\underline{\lambda}}^{(2)}$ defined by:
\begin{itemize}
\item if $\lambda_1=\lambda_2=\lambda$ for some $\lambda>0$,
$$I_{\mathrm{MD},\underline{\lambda}}^{(2)}(x)
:=((\nu/2)^{\nu/(2-\nu)}-(\nu/2)^{2/(2-\nu)})\left(\frac{x^2}{\lambda}\right)^{1/(2-\nu)};$$
\item if $\lambda_1>\lambda_2$,
$$I_{\mathrm{MD},\underline{\lambda}}^{(2)}(x):=\left\{\begin{array}{ll}
(\nu^{\nu/(1-\nu)}-\nu^{1/(1-\nu)})\left(\frac{x}{\lambda_1-\lambda_2}\right)^{1/(1-\nu)}&\ \mbox{if}\ x\geq 0\\
\infty&\ \mbox{if}\ x<0;
\end{array}\right.$$
\item if $\lambda_1<\lambda_2$,
$$I_{\mathrm{MD},\underline{\lambda}}^{(2)}(x):=\left\{\begin{array}{ll}
(\nu^{\nu/(1-\nu)}-\nu^{1/(1-\nu)})\left(\frac{x}{-(\lambda_2-\lambda_1)}\right)^{1/(1-\nu)}&\ \mbox{if}\ x\leq 0\\
\infty&\ \mbox{if}\ x>0.
\end{array}\right.$$
\end{itemize}
\end{proposition}

\begin{remark}\label{rem:weakconv-NCMD-type2}
The sets $\{x\in\mathbb{R}:I_{\mathrm{MD},\underline{\lambda}}^{(2)}(x)<\infty\}$ (see Proposition \ref{prop:NCMD-type2}) 
coincide with the supports of the weak limits in Proposition \ref{prop:weak-convergence-type2}: we mean $\mathbb{R}$ if 
$\lambda_1=\lambda_2$, $[0,\infty)$ if $\lambda_1>\lambda_2$, and $(-\infty,0]$ if $\lambda_1<\lambda_2$.
\end{remark}

\begin{remark}\label{rem:symmetry-type2}
Assume that $\lambda_1=\lambda_2$. Then: the rate function $I_{\mathrm{LD}}^{(2)}$ in Proposition \ref{prop:LD-type2} is a 
symmetric function (we can say this, even if an explicit expression of $I_{\mathrm{LD}}^{(2)}$ is not available, because $\Psi_{\underline{\nu},\underline{\lambda}}^{(2)}$ is a symmetric function); the weak limit in Proposition 
\ref{prop:weak-convergence-type2} is a symmetric random variable; the rate function $I_{\mathrm{MD}}^{(2)}$ in Proposition 
\ref{prop:NCMD-type2} is a symmetric function.
\end{remark}

\section{Comparisons between rate functions}\label{sec:comparisons}
In this section we compare the rate functions for the two types of fractional Skellam process, and for different values
of $\nu$. Moreover, we present some plots.

\subsection{Results and remarks}
All the LDPs presented in the previous sections are governed by rate functions that uniquely vanish at $x=0$. So, as we 
explain below, it is interesting to compare the rate functions, at least around $x=0$. 

\begin{remark}\label{rem:simplified-notation-psi}
Throughout this section we always assume that $\nu_1=\nu_2=\nu$ for some $\nu\in(0,1)$. So we simply write
$\Psi_{\nu,\underline{\lambda}}^{(1)}$ in place of the function $\Psi_{\underline{\nu},\underline{\lambda}}^{(1)}$ in 
Proposition \ref{prop:LD-type1} (see \eqref{eq:LD-GE-limit-type1}).
\end{remark}

We start by comparing $I_{\mathrm{LD}}^{(1)}$ in Proposition \ref{prop:LD-type1} and $I_{\mathrm{LD}}^{(2)}$ in Proposition 
\ref{prop:LD-type2}.

\begin{proposition}\label{prop:comparison-rf-LD}
Assume that $\nu_1=\nu_2=\nu$ for some $\nu\in(0,1)$. Then $I_{\mathrm{LD}}^{(1)}(0)=I_{\mathrm{LD}}^{(2)}(0)=0$ and, for
$x\neq 0$, we have $I_{\mathrm{LD}}^{(2)}(x)>I_{\mathrm{LD}}^{(1)}(x)>0$.
\end{proposition}
\begin{proof}
We start noting that, for the function $\Psi_{\nu,\underline{\lambda}}^{(1)}$ in Proposition \ref{prop:LD-type1} (see Remark 
\ref{rem:simplified-notation-psi} and \eqref{eq:LD-GE-limit-type1}) and the function 
$\Psi_{\nu,\underline{\lambda}}^{(2)}$ in Proposition \ref{prop:LD-type2} (see \eqref{eq:LD-GE-limit-type2} respectively), 
we have
\begin{equation}\label{eq:equivalent-inequality}
\mbox{$\Psi_{\nu,\underline{\lambda}}^{(1)}(0)=\Psi_{\nu,\underline{\lambda}}^{(2)}(0)=0$;
$\Psi_{\nu,\underline{\lambda}}^{(1)}(\theta)>\Psi_{\nu,\underline{\lambda}}^{(2)}(\theta)$ for $\theta\neq 0$.}
\end{equation}
The first statement in \eqref{eq:equivalent-inequality} is immediate. For the second statement we have two cases (for 
completeness we remark that $\lambda_1(e^\theta-1)+\lambda_2(e^{-\theta}-1)\geq 0$ for all $\theta\in\mathbb{R}$ if 
$\lambda_1=\lambda_2$).
\begin{itemize}
\item If $\theta>0$, then
$$\lambda_1(e^\theta-1)>\max\{\lambda_1(e^\theta-1)+\lambda_2(e^{-\theta}-1),0\}$$
which yields $\Psi_{\nu,\underline{\lambda}}^{(1)}(\theta)>\Psi_{\nu,\underline{\lambda}}^{(2)}(\theta)$
by taking the power with exponent $1/\nu$;
\item If $\theta<0$, then
$$\lambda_2(e^{-\theta}-1)>\max\{\lambda_1(e^\theta-1)+\lambda_2(e^{-\theta}-1),0\}$$
which yields $\Psi_{\nu,\underline{\lambda}}^{(1)}(\theta)>\Psi_{\nu,\underline{\lambda}}^{(2)}(\theta)$
again, by taking the power with exponent $1/\nu$.
\end{itemize}
Thus \eqref{eq:equivalent-inequality} is checked.

We also remark that, for every $x\in\mathbb{R}$, there exist $\theta_x^{(1)},\theta_x^{(2)}\in\mathbb{R}$ such that
$$I_{\mathrm{LD}}^{(1)}(x)=\theta_x^{(1)}x-\Psi_{\nu,\underline{\lambda}}^{(1)}(\theta_x^{(1)})
\ \mbox{and}\ I_{\mathrm{LD}}^{(2)}(x)=\theta_x^{(2)}x-\Psi_{\nu,\underline{\lambda}}^{(2)}(\theta_x^{(2)});$$
moreover $\theta_x^{(1)}=\theta_x^{(2)}=0$ if $x=0$, and $\theta_x^{(1)},\theta_x^{(2)}\neq 0$ if $x\neq 0$. Then, if
$x\neq 0$, we get
$$I_{\mathrm{LD}}^{(1)}(x)=\theta_x^{(1)}x-\Psi_{\nu,\underline{\lambda}}^{(1)}(\theta_x^{(1)})
<\theta_x^{(1)}x-\Psi_{\nu,\underline{\lambda}}^{(2)}(\theta_x^{(1)})
\leq\sup_{\theta\in\mathbb{R}}\{\theta x-\Psi_{\nu,\underline{\lambda}}^{(2)}(\theta)\}
=I_{\mathrm{LD}}^{(2)}(x),$$
where the strict inequality holds by \eqref{eq:equivalent-inequality}, and by taking into account that $\theta_x^{(1)}\neq 0$.
\end{proof}

The next Proposition \ref{prop:comparison-rf-MD} provides a similar result which concerns the comparison of $I_{\mathrm{MD}}^{(1)}$
in Proposition \ref{prop:NCMD-type1} and $I_{\mathrm{MD},\underline{\lambda}}^{(2)}$ in Proposition \ref{prop:NCMD-type2}. In 
particular, it is possible to obtain the same strict inequality, for all $x\neq 0$, only if $\lambda_1\neq\lambda_2$ (note that in 
such a case $\alpha_1(\nu)=\alpha_2(\nu)=1-\nu$ by \eqref{eq:exponent} and \eqref{eq:exponents}).

\begin{proposition}\label{prop:comparison-rf-MD}
We have $I_{\mathrm{MD}}^{(1)}(0)=I_{\mathrm{MD},\underline{\lambda}}^{(2)}(0)=0$. Moreover, if $x\neq 0$, we have
two cases.
\begin{enumerate}
\item If $\lambda_1\neq\lambda_2$, then $I_{\mathrm{MD},\underline{\lambda}}^{(2)}(x)>I_{\mathrm{MD}}^{(1)}(x)>0$.
\item If $\lambda_1=\lambda_2=\lambda$ for some $\lambda>0$, there exists $\delta_{\nu,\lambda}>0$ such that:
$I_{\mathrm{MD},\underline{\lambda}}^{(2)}(x)>I_{\mathrm{MD}}^{(1)}(x)>0$ if $0<|x|<\delta_{\nu,\lambda}$,
$I_{\mathrm{MD}}^{(1)}(x)>I_{\mathrm{MD},\underline{\lambda}}^{(2)}(x)>0$ if $|x|>\delta_{\nu,\lambda}$,
and $I_{\mathrm{MD},\underline{\lambda}}^{(2)}(x)=I_{\mathrm{MD}}^{(1)}(x)>0$ if $|x|=\delta_{\nu,\lambda}$.
\end{enumerate}
\end{proposition}
\begin{proof}
The equalities $I_{\mathrm{MD}}^{(1)}(0)=I_{\mathrm{MD},\underline{\lambda}}^{(2)}(0)=0$ (case $x=0$) are immediate.
So, in what follows, we take $x\neq 0$. We start with the case $\lambda_1\neq\lambda_2$, and we have two cases.
\begin{itemize}
\item Assume that $\lambda_1>\lambda_2$. Then for $x<0$ we have $I_{\mathrm{MD}}^{(1)}(x)<\infty=I_{\mathrm{MD},\underline{\lambda}}^{(2)}(x)$.
For $x>0$ we have $\frac{x}{\lambda_1}<\frac{x}{\lambda_1-\lambda_2}$, which is trivially equivalent to
$I_{\mathrm{MD}}^{(1)}(x)<I_{\mathrm{MD},\underline{\lambda}}^{(2)}(x)$.
\item Assume that $\lambda_1<\lambda_2$. Then for $x>0$ we have $I_{\mathrm{MD}}^{(1)}(x)<\infty=I_{\mathrm{MD},\underline{\lambda}}^{(2)}(x)$.
For $x<0$ we have $\frac{x}{-\lambda_2}<\frac{x}{-(\lambda_2-\lambda_1)}$, which is trivially equivalent to
$I_{\mathrm{MD}}^{(1)}(x)<I_{\mathrm{MD},\underline{\lambda}}^{(2)}(x)$.
\end{itemize}
Finally, if $\lambda_1=\lambda_2=\lambda$ for some $\lambda>0$, the statement to prove trivially holds noting that, 
for two constants $c_{\nu,\lambda}^{(1)},c_{\nu,\lambda}^{(2)}>0$, we have
$I_{\mathrm{MD}}^{(1)}(x)=c_{\nu,\lambda}^{(1)}|x|^{1/(1-\nu)}$ and 
$I_{\mathrm{MD},\underline{\lambda}}^{(2)}(x)=c_{\nu,\lambda}^{(2)}|x|^{1/(1-\nu/2)}$.
\end{proof}

Proposition \ref{prop:comparison-rf-LD} tells us that if we compare the rate functions in Propositions \ref{prop:LD-type1} 
and \ref{prop:LD-type2}, the rate function of the fractional Skellam process of type 2 is larger than that of the 
fractional Skellam process of type 1. Proposition \ref{prop:comparison-rf-MD} tells us the same for the rate functions in 
Propositions \ref{prop:NCMD-type1} and \ref{prop:NCMD-type2} but, when $\lambda_1=\lambda_2$, the rate function of the fractional 
Skellam process of type 2 is larger only around $x=0$.

These inequalities between rate functions allow us to say that the convergence of random variables for the fractional 
Skellam process of type 2 is faster than the corresponding convergence for the fractional Skellam process of type 1. We explain 
this by considering the LDPs in Proposition \ref{prop:LD-type1}, with $\nu_1=\nu_2=\nu$ for some $\nu\in(0,1)$, and in Proposition 
\ref{prop:LD-type2}. Indeed, for every $\delta>0$ we have
$$\lim_{t\to\infty}\frac{1}{t}\log P\left(\frac{|Y_{\nu,\underline{\lambda}}(t)|}{t}>\delta\right)=-J_{\mathrm{LD}}^{(1)}(\delta),
\ \mbox{where}\ J_{\mathrm{LD}}^{(1)}(\delta):=\min\{I_{\mathrm{LD}}^{(1)}(\delta),I_{\mathrm{LD}}^{(1)}(-\delta)\}$$
and
$$\lim_{t\to\infty}\frac{1}{t}\log P\left(\frac{|Z_{\nu,\underline{\lambda}}(t)|}{t}>\delta\right)=-J_{\mathrm{LD}}^{(2)}(\delta),
\ \mbox{where}\ J_{\mathrm{LD}}^{(2)}(\delta):=\min\{I_{\mathrm{LD}}^{(2)}(\delta),I_{\mathrm{LD}}^{(2)}(-\delta)\};$$
therefore $J_{\mathrm{LD}}^{(2)}(\delta)>J_{\mathrm{LD}}^{(1)}(\delta)>0$ and, for every $\varepsilon\in \left(0,J_{\mathrm{LD}}^{(2)}(\delta)-J_{\mathrm{LD}}^{(1)}(\delta)\right)$, there exists $t_\varepsilon$ such that
$$\frac{P\left(\frac{|Z_{\nu,\underline{\lambda}}(t)|}{t}>\delta\right)}{P\left(\frac{|Y_{\nu,\underline{\lambda}}(t)|}{t}>\delta\right)}
<e^{-t(J_{\mathrm{LD}}^{(2)}(\delta)-J_{\mathrm{LD}}^{(1)}(\delta)-\varepsilon)}\ \mbox{for all}\ t>t_\varepsilon,$$
where $e^{-t(J_{\mathrm{LD}}^{(2)}(\delta)-J_{\mathrm{LD}}^{(1)}(\delta)-\varepsilon)}\to 0$ as $t\to\infty$.

We can follow the same lines to obtain a similar estimate starting from LDPs in Propositions \ref{prop:NCMD-type1}
and \ref{prop:NCMD-type2}; in this case, when $\lambda_1=\lambda_2$, $\delta$ has to be chosen small enough (see Proposition 
\ref{prop:comparison-rf-MD}). Here we do not repeat all the computations; however, we can say that if we set
$$J_{\mathrm{MD}}^{(1)}(\delta):=\min\{I_{\mathrm{MD}}^{(1)}(\delta),I_{\mathrm{MD}}^{(1)}(-\delta)\}
\ \mbox{and}\ J_{\mathrm{MD}}^{(2)}(\delta):=\min\{I_{\mathrm{MD},\underline{\lambda}}^{(2)}(\delta),
I_{\mathrm{MD},\underline{\lambda}}^{(2)}(-\delta)\},$$
we have $J_{\mathrm{MD}}^{(2)}(\delta)>J_{\mathrm{MD}}^{(1)}(\delta)>0$ and, for every $\varepsilon\in \left(0,J_{\mathrm{MD}}^{(2)}(\delta)-J_{\mathrm{MD}}^{(1)}(\delta)\right)$, there exists $t_\varepsilon$ such that
$$\frac{P\left(\frac{(a_tt)^{\alpha_2(\nu)}|Z_{\nu,\underline{\lambda}}(t)|}{t}>\delta\right)}
{P\left(\frac{(a_tt)^{\alpha_1(\nu)}|Y_{\nu,\underline{\lambda}}(t)|}{t}>\delta\right)}
<e^{-(J_{\mathrm{MD}}^{(2)}(\delta)-J_{\mathrm{MD}}^{(1)}(\delta)-\varepsilon)/a_t}\ \mbox{for all}\ t>t_\varepsilon,$$
where $e^{-(J_{\mathrm{MD}}^{(2)}(\delta)-J_{\mathrm{MD}}^{(1)}(\delta)-\varepsilon)/a_t}\to 0$ as $t\to\infty$.

\begin{remark}\label{rem:connections-with-variances-in-KerssLeonenkoSikorskii}
We can say that it is not surprising that the convergence of random variables for the fractional Skellam process of type 2 is 
faster than the corresponding convergence for the fractional Skellam process of type 1. Indeed the fractional Skellam process of
type 1 is defined with two independent random time-changes for the involved fractional Poisson processes; on the contrary,
for the fractional Skellam process of type 2, we have a unique independent random time-change. So, in some sense, in the second
case we have less randomness. Moreover, if we compare the asymptotic behavior (as $t\to\infty$) of the variances provided by 
Remarks 3.1 and 3.2 in \cite{KerssLeonenkoSikorskii}, i.e.
$$\mathrm{Var}[Y_{\nu,\underline{\lambda}}(t)]=\frac{t^\nu(\lambda_1+\lambda_2)}{\Gamma(1+\nu)}+\frac{(\lambda_1^2+\lambda_2^2)t^{2\nu}}{\nu}
\left(\frac{1}{\Gamma(2\nu)}-\frac{1}{\nu\Gamma^2(\nu)}\right)$$
and
$$\mathrm{Var}[Z_{\nu,\underline{\lambda}}(t)]=\frac{t^\nu(\lambda_1+\lambda_2)}{\Gamma(1+\nu)}+(\lambda_1-\lambda_2)^2t^{2\nu}
\left(\frac{2}{\Gamma(2\nu+1)}-\frac{1}{\Gamma^2(1+\nu)}\right),$$
in the second case we have smaller asymptotic variances. Indeed we have
$$\lim_{t\to\infty}\frac{\mathrm{Var}[Y_{\nu,\underline{\lambda}}(t)]}{t^{2\nu}}>
\lim_{t\to\infty}\frac{\mathrm{Var}[Z_{\nu,\underline{\lambda}}(t)]}{t^{2\nu}},$$
noting that
$$\lim_{t\to\infty}\frac{\mathrm{Var}[Y_{\nu,\underline{\lambda}}(t)]}{t^{2\nu}}=
\frac{\lambda_1^2+\lambda_2^2}{\nu}\left(\frac{1}{\Gamma(2\nu)}-\frac{1}{\nu\Gamma^2(\nu)}\right)$$
and
$$\lim_{t\to\infty}\frac{\mathrm{Var}[Z_{\nu,\underline{\lambda}}(t)]}{t^{2\nu}}=
(\lambda_1-\lambda_2)^2\left(\frac{2}{\Gamma(2\nu+1)}-\frac{1}{\Gamma^2(1+\nu)}\right),$$
where
$$\frac{2}{\Gamma(2\nu+1)}-\frac{1}{\Gamma^2(1+\nu)}
=\frac{2}{2\nu\Gamma(2\nu)}-\frac{1}{\nu^2\Gamma^2(\nu)}
=\frac{1}{\nu}\left(\frac{1}{\Gamma(2\nu)}-\frac{1}{\nu\Gamma^2(\nu)}\right)$$
and $(\lambda_1-\lambda_2)^2<\lambda_1^2+\lambda_2^2$.
\end{remark}

Now we consider comparisons between rate functions for different values of $\nu\in(0,1)$.
We mean the rate functions in Propositions \ref{prop:LD-type1} and \ref{prop:LD-type2}, and
we restrict our attention to a comparison around $x=0$. In view of what follows, we consider
some slightly different notation: $I_{\mathrm{LD},\nu}^{(1)}$ in place of $I_{\mathrm{LD}}^{(1)}$
in Proposition \ref{prop:LD-type1}, with $\nu_1=\nu_2=\nu$ for some $\nu\in(0,1)$, and 
$I_{\mathrm{LD},\nu}^{(2)}$ in place of $I_{\mathrm{LD}}^{(2)}$ in Proposition \ref{prop:LD-type2}.

\begin{proposition}\label{prop:comparison-rf-LD-nu}
Let $\nu,\eta\in(0,1)$ be such that $\eta<\nu$. Then, for $k\in\{1,2\}$, we have 
$I_{\mathrm{LD},\eta}^{(k)}(0)=I_{\mathrm{LD},\nu}^{(k)}(0)=0$ and, for some $\delta>0$,
$I_{\mathrm{LD},\eta}^{(k)}(x)>I_{\mathrm{LD},\nu}^{(k)}(x)>0$ for $0<|x|<\delta$.
\end{proposition}
\begin{proof}
We take an arbitrarily fixed $k\in\{1,2\}$. Then, for every $x\in\mathbb{R}$, there exists
$\theta_x^{(\nu,k)}\in\mathbb{R}$ such that
$$I_{\mathrm{LD},\nu}^{(k)}(x)=\theta_x^{(\nu,k)}x-\Psi_{\nu,\underline{\lambda}}^{(k)}(\theta_x^{(\nu,k)});$$
we recall that we are referring to the function $\Psi_{\nu,\underline{\lambda}}^{(1)}$ in 
\eqref{eq:LD-GE-limit-type1} (see also Remark \ref{rem:simplified-notation-psi}) and to the function 
$\Psi_{\nu,\underline{\lambda}}^{(2)}$ in \eqref{eq:LD-GE-limit-type2}. We can also say that $\theta_x^{(\nu,k)}=0$
if and only if $x=0$ and, moreover, there exists $\delta>0$ such that 
$0\leq\Psi_{\nu,\underline{\lambda}}^{(1)}(\theta_x^{(\nu,k)})<1$ if $|x|<\delta$. Then, by taking into account
the same formulas with $\eta$ in place of $\nu$, it is easy to check that
$$0\leq\Psi_{\eta,\underline{\lambda}}^{(k)}(\theta_x^{(\nu,k)})<\Psi_{\nu,\underline{\lambda}}^{(k)}(\theta_x^{(\nu,k)})<1$$
(see \eqref{eq:LD-GE-limit-type1} and \eqref{eq:LD-GE-limit-type2}; moreover we take into account that 
$\frac{1}{\eta}>\frac{1}{\nu}$), whence we obtain
$$I_{\mathrm{LD},\nu}^{(k)}(x)=\theta_x^{(\nu,k)}x-\Psi_{\nu,\underline{\lambda}}^{(k)}(\theta_x^{(\nu,k)})
<\theta_x^{(\nu,k)}x-\Psi_{\eta,\underline{\lambda}}^{(k)}(\theta_x^{(\nu,k)})
\leq\sup_{\theta\in\mathbb{R}}\{\theta x-\Psi_{\eta,\underline{\lambda}}^{(k)}(\theta)\}=I_{\mathrm{LD},\eta}^{(k)}(x).$$
This completes the proof.
\end{proof}

As a consequence of Proposition \ref{prop:comparison-rf-LD-nu} we can present some estimates that allow us to compare
the convergence to zero of different families of random variables for different values of $\nu\in(0,1)$. Roughly speaking 
we can say the smaller the $\nu$, the faster the convergence of the random variables to zero. This could be explained by
presenting a modified version of the computations presented just after the proof of Proposition \ref{prop:comparison-rf-MD};
here we omit the details.

\subsection{Some plots}
\begin{figure}[!ht]
\centering
\includegraphics[width=\linewidth]{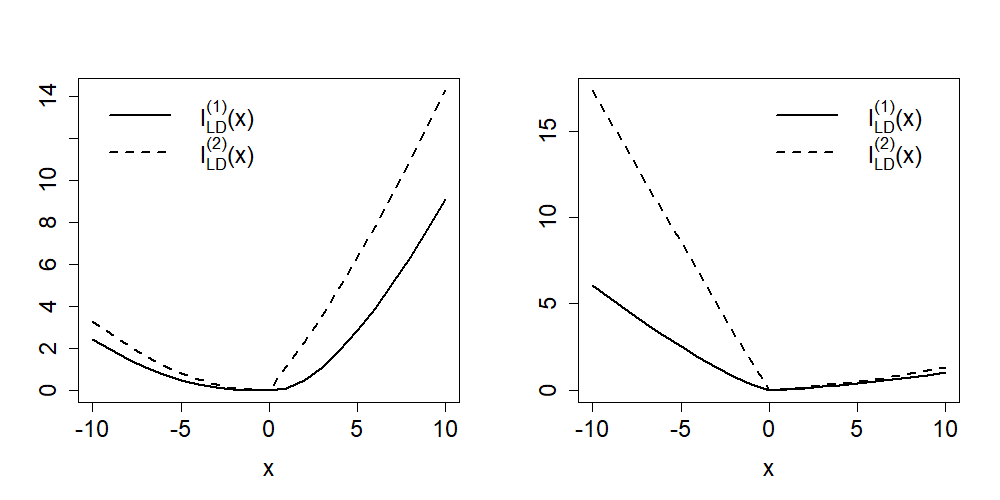}
\includegraphics[width=\linewidth]{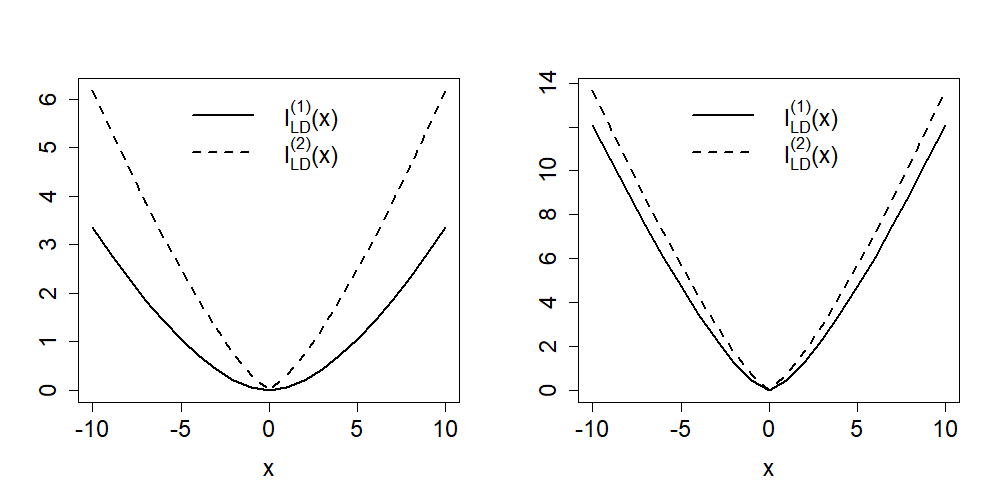}
\caption{Top left:  $(\lambda_1,\lambda_2,\nu)=(1,3,.7)$. Top right:  $(\lambda_1,\lambda_2,\nu)=(5,1,.3)$. 
Bottom left: $(\lambda_1,\lambda_2,\nu)=(2,2,.5)$. Bottom right:  $(\lambda_1,\lambda_2,\nu)=(.5,.5,.5)$.}
\label{fig1}
\end{figure}

We start with Figure \ref{fig1} which shows that the rate functions $I_{\mathrm{LD}}^{(2)}(x)$ and $I_{\mathrm{LD}}^{(1)}(x)$ 
when $\nu_1=\nu_2=\nu$ for various sets of parameters $(\lambda_1,\lambda_2,\nu)$. In each case the plots agree with the 
statements in Proposition \ref{prop:comparison-rf-LD}, i.e. $I_{\mathrm{LD}}^{(2)}(x)>I_{\mathrm{LD}}^{(1)}(x)>0$ for $x\neq 0$
and $I_{\mathrm{LD}}^{(2)}(0)=I_{\mathrm{LD}}^{(1)}(0)=0$. In particular, if $\lambda_1=\lambda_2$, the rate functions are 
symmetric (around zero) as announced in Remarks \ref{rem:symmetry-type1} and \ref{rem:symmetry-type2}.

\begin{figure}[!ht]
\centering
\includegraphics[width=\linewidth]{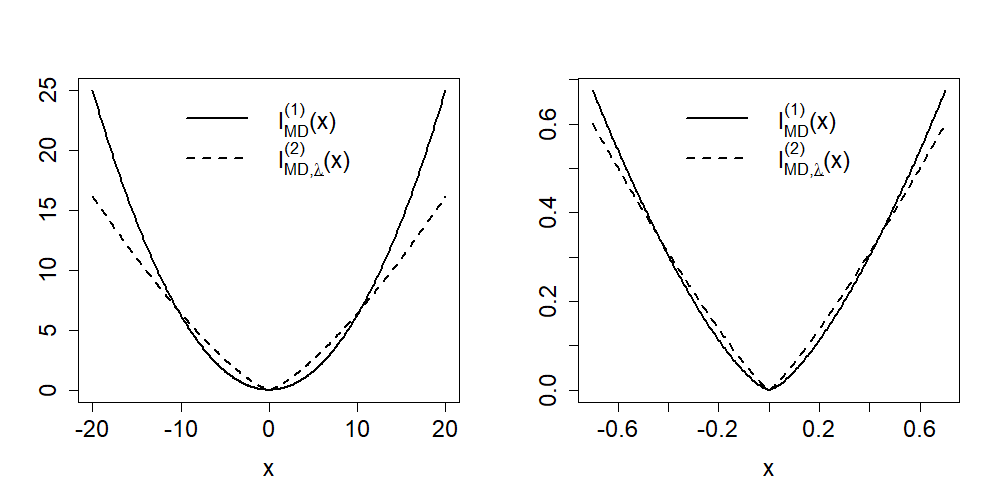}
\caption{Left: $\lambda_1=\lambda_2=2, \nu=.5$. Right:  $\lambda_1=\lambda_2=.5, \nu=.3$.}
\label{fig2}
\end{figure}

In Figure \ref{fig2} we present the plots of $I_{\mathrm{MD}}^{(1)}(x)$ and $I_{\mathrm{MD},\underline{\lambda}}^{(2)}(x)$
when $\lambda_1=\lambda_2$. For each set of parameters, the plots agree with a statement of Proposition \ref{prop:comparison-rf-MD}:
$I_{\mathrm{MD}}^{(1)}(x)$ is smaller than $I_{\mathrm{MD},\underline{\lambda}}^{(2)}(x)$ when $x$ is within a bounded interval that
includes zero, and becomes larger outside of that interval.

\begin{figure}[!ht]
\centering
\includegraphics[width=\linewidth]{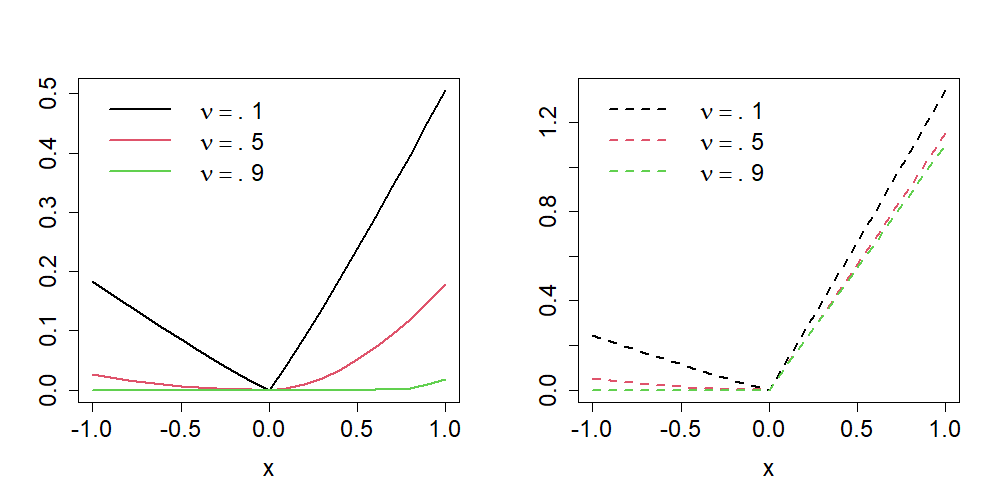}
\caption{Left: $I_{\mathrm{LD}}^{(1)}(x)$ with $\lambda_1=1, \lambda_2=3$. Right: $I_{\mathrm{LD}}^{(2)}(x)$ with $\lambda_1=1, \lambda_2=3$.}
\label{fig3}
\end{figure}

Finally in Figure \ref{fig3} we present the plots of $I_{\mathrm{LD}}^{(1)}(x)$ and $I_{\mathrm{LD}}^{(2)}(x)$ for
different values of $\nu$. These plots agree with the statements in Proposition \ref{prop:comparison-rf-LD-nu}.

\section{Concluding remarks}\label{sec:concluding-remarks}
In this paper we prove noncentral moderate deviations for two fractional Skellam processes presented in the literature.
The main tool used in the proofs is the G\"artner Ellis Theorem (see Theorem \ref{th:GE}) which can be applied because
the involved moment generating functions are available in a closed form given in terms of the Mittag-Leffler function.

For the classical (non-fractional) Skellam process we can obtain a classical moderate deviation result that fills the
gap between the two following asymptotic regimes:
\begin{enumerate}
	\item the convergence of $\frac{S_{\underline{\lambda}}(t)}{t}-(\lambda_1+\lambda_2)$ in probability to zero, which
	is governed by a LDP with speed $t$; 
	\item the weak convergence of $\sqrt{t}\left(\frac{S_{\underline{\lambda}}(t)}{t}-(\lambda_1+\lambda_2)\right)$ to
	a centered Normal distribution.
\end{enumerate}

In general the applications of the G\"artner Ellis Theorem for random time-changed processes are quite standard. However, 
in order to obtain noncentral moderate deviation results as the ones in this paper, it is important to have a
random time-change with a slowing effect as happens in this paper with the inverse of the stable subordinator; in this
way we have normalized processes that tend to zero (as happens for $\frac{N_{\nu_1,\lambda_1}(t)-N_{\nu_2,\lambda_2}(t)}{t}$
and $\frac{S_{\underline{\lambda}}(L_\nu(t))}{t}$ in this paper).
In future work one could consider more general models with (independent) random time-changes in terms of an inverse
of a general subordinator (see e.g. the recent reference \cite{GuptaMaheshwari}) with the same slowing effect.

The results in this work (and in \cite{BeghinMacciSPL2022}) concern one-dimensional cases. It could be nice to obtain 
sample-path versions of these results. One could try to combine the sample-path results for light-tailed Lévy processes
in \cite{deAcosta} (which can be applied to non-fractional Skellam processes) with possible sample-path results
for the inverse of the stable subordinator; unfortunately the derivation of sample-path results for the inverse of the 
stable subordinator seems to be a difficult task.

We conclude with a discussion on the connection between the scaling exponents of the normalizing factors for the moderate 
deviations of the fractional Skellam processes, and the concept of \emph{long-range dependence} (LRD). We recall (see e.g. 
Definition 2.1 in \cite{KatariaKhandakarJTP2021}) that a non-stationary stochastic process $\{X(t):t\geq 0\}$ has the 
long-range dependence if we have the following property for the correlation coefficient:
$$\frac{\mathrm{Cov}(X(t),X(s))}{\sqrt{\mathrm{Var}[X(t)]\mathrm{Var}[X(s)]}}\sim c(s)t^{-h}\quad (\mbox{as}\ t\to\infty),$$
for some $c(s)>0$ and $h\in(0,1)$. Then we have the following statements on the scaling exponents  $\alpha_1(\nu)$ in 
\eqref{eq:exponent} and  $\alpha_2(\nu)$ in \eqref{eq:exponents}.
\begin{itemize}
	\item For the Fractional Skellam process of type 1 with $\nu_1=\nu_2=\nu$ for some $\nu\in(0,1)$, we have the LRD for $\{Y_{\nu,\underline{\lambda}}(t):t\geq 0\}$ with $h=1-\alpha_1(\nu)=\nu$ (this is a consequence of eqs. (19) and (20) and
	Remark 2 in \cite{KatariaKhandakarJTP2021}, and some computations).
	\item For the Fractional Skellam process of type 2 we have the LRD for $\{Z_{\nu,\underline{\lambda}}(t):t\geq 0\}$ with
	$$h=1-\alpha_2(\nu)=\left\{\begin{array}{ll}
		\nu/2&\ \mbox{if}\ \lambda_1=\lambda_2\\
		\nu&\ \mbox{if}\ \lambda_1\neq\lambda_2
	\end{array}\right.$$
	(this is a consequence of eqs. (2.11)-(2.16) in \cite{KatariaKhandakar-arxiv} with $k=1$, and some computations).
\end{itemize}

\paragraph{Funding.}
CM acknowledges the support of MIUR Excellence Department Project awarded to the Department of Mathematics, University
of Rome Tor Vergata (CUP E83C18000100006 and CUP E83C23000330006), of University of Rome Tor Vergata (project 
"Asymptotic Methods in Probability" (CUP E89C20000680005) and project "Asymptotic Properties in Probability" (CUP 
E83C22001780005)) and of Indam-GNAMPA.

\paragraph{Acknowledgements.} 
The authors thank two referees for their positive opinions, and for some useful comments.


\begin{thebibliography}{spc}
\bibitem{BeghinMacciSPL2013}
L. Beghin, C. Macci (2013) Large deviations for fractional Poisson processes. Statist.
Probab. Lett. 83, 1193--1202.
\bibitem{BeghinMacciSPL2017}
L. Beghin, C. Macci (2017) Asymptotic results for a multivariate version of the alternative
fractional Poisson process. Statist. Probab. Lett. 129, 260--268.
\bibitem{BeghinMacciSPL2022}
L. Beghin, C. Macci (2022) Non-central moderate deviations for compound fractional Poisson 
processes. Statist. Probab. Lett. 185, Paper No. 109424, 8 pp.
\bibitem{BeghinMacciMartinucci}
L. Beghin, C. Macci, B. Martinucci (2021) Random time-changes and asymptotic results for a 
class of continuous-time Markov chains on integers with alternating rates. Mod. Stoch. Theory
Appl. 8, 63--91.
\bibitem{BeghinOrsingher2009}
L. Beghin, E. Orsingher (2009) Fractional Poisson processes and related planar random motions. 
Electron. J. Probab. 14, 1790--1827.
\bibitem{BeghinOrsingher2010}
L. Beghin, E. Orsingher (2010) Poisson-type processes governed by fractional and higher-order 
recursive differential equations. Electron. J. Probab. 15, 684--709.
\bibitem{BuchakSakhno2017}
K. Buchak, L. Sakhno (2017) Compositions of Poisson and gamma processes. Mod. Stoch. Theory 
Appl. 4, 161--188.
\bibitem{BuchakSakhno2019}
K.V. Buchak, L.M. Sakhno (2019) On the governing equations for Poisson and Skellam processes 
time-changed by inverse subordinators. Theory Probab. Math. Statist. 98, 91--104.
\bibitem{deAcosta}
A. de Acosta (1994) Large deviations for vector-valued Lévy processes. Stochastic Process.
Appl. 51, 75--115.
\bibitem{DemboZeitouni}
A. Dembo, O. Zeitouni (1998) Large Deviations Techniques and Applications, 2nd  edn. Springer.
\bibitem{GarraOrsingherPolito}
R. Garra,  E. Orsingher, F. Polito (2015) State dependent fractional point processes. J. Appl.
Probab. 52, 18--36 (2015).
\bibitem{GiulianoMacci}
R. Giuliano, C. Macci (2023) Some examples of noncentral moderate deviations for sequences of
real random variables. Mod. Stoch. Theory Appl. 10, 111--144.
\bibitem{GorenfloKilbasMainardiRogosin}
R. Gorenflo, A.A. Kilbas, F. Mainardi, S.V. Rogosin (2014) Mittag-Leffler Functions, Related
Topics and Applications. Springer, New York.
\bibitem{GuptaKumarLeonenko}
N. Gupta, A. Kumar, N. Leonenko (2020) Fractional Skellam process of order $k$ and beyond.
Entropy 22, Paper No. 1193.
\bibitem{GuptaMaheshwari}
N. Gupta, A. Maheshwari (2023+) Fractional generalizations of the compound Poisson processes.
Available at \texttt{https://arxiv.org/pdf/2307.12252.pdf}
\bibitem{KatariaKhandakarJTP2021}
K.K. Kataria, M. Khandakar (2021) On the long-range dependence of mixed fractional Poisson process.
J. Theoret. Probab. 34, 1607--1622. 
\bibitem{KatariaKhandakar-arxiv}
K.K. Kataria, M. Khandakar (2021+) Fractional Skellam process of order $k$. Available at
\texttt{https://arxiv.org/pdf/2103.09187.pdf}
\bibitem{KatariaKhandakarFCAA2022}
K.K. Kataria, M. Khandakar (2022) Skellam and time-changed variants of the generalized fractional 
counting processes. Fract. Calc. Appl. Anal. 25, 1873--1907.
\bibitem{KatariaKhandakarSAA2022}
K.K. Kataria, M. Khandakar (2022) Time-changed space-time fractional Poisson process. 
Stoch. Anal. Appl. 40, 246--267.
\bibitem{KerssLeonenkoSikorskii}
A. Kerss, N.N. Leonenko, A. Sikorskii (2014) Fractional Skellam processes with applications 
to finance. Fract. Calc. Appl. Anal. 17, 532--551.
\bibitem{MacciPacchiarottiVilla}
C. Macci, B. Pacchiarotti, E.Villa (2022) Asymptotic results for families of random variables having
power series distributions. Mod. Stoch. Theory Appl. 9, 207--228. 
\bibitem{MeerschaertNaneVellaisamy}
M.M. Meerschaert, E. Nane, P. Vellaisamy (2011) The fractional Poisson process and the inverse 
stable subordinator. Electron. J. Probab. 16, 1600--1620. 
\bibitem{OrsingherPolito}
E. Orsingher, F. Polito (2012) The space-fractional Poisson process. Statist.
Probab. Lett. 82, 852--858.
\bibitem{PoganyTomovski}
T.K. Pogány, Ž. Tomovski (2016) Probability distribution built by Prabhakar function. Related
Turán and Laguerre inequalities. Integral Transforms Spec. Funct. 27, 783--793.
\end{thebibliography}
\end{document}